\documentclass[11pt]{article}
\pdfoutput=1

\usepackage{graphicx}
\usepackage[english]{babel}
\usepackage{hyperref}
\usepackage{amsmath}

\usepackage{tocloft}

\title{Bernstein- and Markov-type inequalities}

\author{Sergei Kalmykov, B\'ela Nagy and Vilmos Totik}

\date{May 19, 2021}

\def\abstracttext{This survey discusses the classical
Bernstein and Markov
inequalities for the derivatives of polynomials, as well as some
of their extensions to general sets.}

\def\MSCnumbers{26D05, 42A05} 

\newtheorem{thh}{Theorem}

\newtheorem{problem}[thh]{Problem}

\def\F{\Phi}
\def\G{\Gamma}

\newcommand{\R}{\mathbf{R}}
\newcommand{\C}{\mathbf{C}}

\newcommand{\ee}{\mathrm{e}}
\newcommand{\ii}{\mathrm{i}}

\begin{document}

\maketitle

\begin{abstract}
\abstracttext

MSC Classification: \MSCnumbers
\end{abstract}

\tableofcontents

\section{The original Bernstein and Markov inequalities}\label{sect1}
In 1912 S. N. Bernstein proved in \cite{Bernstein}
his famous inequality
that now takes the form
\begin{equation}
|T_n'(\theta)|\le n\sup_t|T_n(t)|,\qquad \theta\in\R,
\label{b0}
\end{equation}
where
\[
T_n(t)=
a_0+(a_1\cos t+b_1\sin t)+\cdots+(a_n\cos nt+b_n\sin nt)
\]
is an arbitrary trigonometric polynomial of degree at most $n$.
In (\ref{b0}) equality occurs for example for $T_n(t)=\sin nt$
and $\theta=0$.
Bernstein stated and proved his inequality in this form
only for even or odd
trigonometric polynomials, and by decomposition into even and odd
parts, for arbitrary trigonometric polynomials
he had $2n$ on the right.
However the
improved version with the correct factor was soon
found by E. Landau and M. Riesz \cite{Riesz},
and L. Fej\'er observed
that the odd case actually implies (\ref{b0})
in its full generality.\footnote{Indeed, it is sufficient
to consider $\theta=0$,
and if we apply Bernstein's theorem for the
odd polynomial
$(T_n(x)-T_n(-x))/2$ at $\theta=0$, then we get (\ref{b0}).}

Let us  rewrite (\ref{b0}) in the form
\[\|T_n'\|\le n\|T_n\|,\]
where $\|T_n\|:=\sup_t|T_n(t)|$ is the supremum norm over
the whole real line.
In general, the supremum norm on a set $E$ is defined as
\[\|f\|_E:=\sup_{t\in E}|f(t)|.\]

If
\[P_n(x)=a_nx^n+a_{n-1}x^{n-1}+\cdots+a_0\]
is an algebraic polynomial of degree at most $n=1,2,\ldots$,
then $P_n(\cos t)$
is a trigonometric polynomial of degree at most $n$, and
for it we obtain from (\ref{b0})
\begin{equation}
|P_n'(x)|\le \frac{n}{\sqrt {1-x^2}} \|P_n\|_{[-1,1]},
\qquad x\in (-1,1),\label{b1}
\end{equation}
which is the ``Bernstein's inequality" for algebraic polynomials.

The right-hand side blows up
as $x\to \pm 1$,  so (\ref{b0}) does not give information on
how large the norm of $P_n'$ can be in terms of the norm of $P_n$.
This question was answered by
the following estimate due
to A. A. Markov \cite{Markov1}
from 1890:
\begin{equation} \|P_n'\|_{[-1,1]}\le n^2\|P_n\|_{[-1,1]}.
\label{m1}\end{equation}

The polynomial inequalities we have been discussing have
various applications.
In approximation
theory they are fundamental in establishing converse results, i.e.,
 when one deduces smoothness from a given rate of
approximation.
For their applications in other areas see \cite{TotikSoz}.

The inequalities (\ref{b1}) and (\ref{m1}) are sharp
and can be applied on any interval
instead of $[-1,1]$.
On more general sets they also give some information,
e.g., if $E=\cup[a_i,b_i]$ consists of finitely
many intervals, then (\ref{m1}) yields
(by applying (\ref{m1}) to each subinterval separately) that
\[
\|P_n'\|_E
\le
n^22\bigl( \min_i (b_i-a_i)\bigr)^{-1}\|P_n\|_E,
\]
but here the ``Markov factor" $2( \min_i (b_i-a_i)\bigr)^{-1}$
on the right is not precise,
it can be replaced by a smaller quantity.

In this paper we shall be interested in the form of
the Bernstein and Markov inequalities on general sets $E$.
The primary concern will be to identify
the best (or asymptotically best) Bernstein and Markov factors
which are connected with geometric
(more precisely, potential theoretic)
properties of the underlying set $E$.
In this respect we mention that until ca. 2000 the analogue
of (\ref{b1}) or (\ref{m1}) was known only in a few special cases,
e.g., for two intervals of equal length,
which can be reduced to the single interval case by the $x\to x^2$
substitution (see \cite{Bor}).
We shall focus on the supremum norm---analogous results
in other norms are scarce, but we shall
mention one in the last section.
Some open problems will also be stated.

There are some interesting local variants of
the Bernstein inequality by V. Andrievskii \cite{A3}
as well as their connection with Bernstein's
and Vasiliev's theorem on approximation of
$|x|$ by polynomials
(\cite{A1}, \cite{A2}, \cite{TotikA}),
 but we shall not discuss them for they need the
concept of Green's functions which we want to avoid in this note.

\section{Equilibrium measures}
Many extensions and generalizations of the original Bernstein and
Markov inequalities
have been found in the last 130 years. We mention here
only one, namely in 1960 V.~S.~ Videnskii \cite{Videnskii1} proved
the analogue of (\ref{b0}) on intervals
shorter than the whole period:  if $\beta\in (0,\pi)$,
then for $\theta\in (-\beta,\beta)$
we have
\begin{equation} |T_n'(\theta)|\le
n\frac{\cos \theta/2}{\sqrt{\sin^{2}\beta/2-\sin^{2}\theta/2}}
\|T_n\|_{[-\beta,\beta]}.
\label{ineq110}\end{equation}
This inequality of Videnskii was sort of a curiosity for
almost half of a century because
the nature of the factor on the right was hidden---we
shall see that it comes from an equilibrium density.
Until recently it was unknown  what the analogue of
the classical inequalities on two (or more intervals),
and even less
on general sets, are, and we shall see that the correct
forms are related to some equilibrium densities.

To formulate the appropriate statements, we need to
introduce a few notions from   potential theory.
For a general reference to logarithmic
potential theory\label{sectpot} see \cite{Ransford}.

Let $E\subset\C$ be a compact subset of the plane.
Think of $E$ as a conductor,
and put a unit charge on $E$, which can freely move in $E$.
After a while the charge
settles, it reaches an equilibrium state where
its internal energy is minimal.
The mathematical formulation is
the following (on the plane, Coulomb's law takes the form that the
repelling force between
charged particles is proportional to the reciprocal of the
distance): except for pathological
cases, there is
a unique probability measure $\mu_E$ on $E$, called the
equilibrium measure
of $E$, that minimizes the energy integral
\begin{equation} \int\int
\log{\frac{1}{|z-t|}}d\mu(z)d\mu(t).
\label{energy}\end{equation}
This $\mu_E$ certainly exists in all the cases we are considering
in this paper.

When $E\subset \R$ we shall denote by $\omega_E(t)$
the density (Radon-Nykodim derivative)
of $\mu_E$ with respect to Lebesgue
measure wherever it exists. It certainly exists in the
(one dimensional) interior of $E$.
For example,
\[ 
\omega_{[-1,1]}(t)=
{\frac{1}{\pi\sqrt{1-t^2}}}, \qquad
t\in(-1,1),\]
is just the well-known Chebyshev distribution.
More generally, if  $E$ consists of finitely many
intervals on the real line, say
\[
E=
\bigcup_{i=1}^m[a_{2i-1},a_{2i}],
\qquad a_1<a_2<\cdots<a_{2m},
\]
then (see \cite{TotikActa})
\begin{equation} \omega_{E}(t)=
\frac{1}{\pi}
\frac{\prod_{i=1}^{m-1}|t -\xi_i|}
{\sqrt{\prod_{i=1}^{2m}{|t -a_i|}}},
\qquad t\in E,
\label{omegaform}\end{equation}
where the $\xi_j\in (a_{2j},a_{2j+1})$,
$j=1,\ldots,m-1$, are the unique solutions of the
system of equations
\begin{equation}
\int_{a_{2j}}^{a_{2j+1}}\frac{\prod_{i=1}^{m-1}(u-\xi_i)}
{\sqrt{\prod_{i=1}^{2m}{|u-a_i|}}}\,du=0,\quad j=1,\ldots,m-1.
\label{7}\end{equation}

In a similar fashion, if $E$ consists of disjoint
smooth Jordan curves and arcs
with arc measure $s_E$, then we set $d\mu_E:=\omega_Eds_E$,
and this $\omega_E$ is then called
the equilibrium
density on $E$. For example, if $E$ is a circle of radius $r$
then $\omega_E(z)\equiv 1/(2\pi r)$ on $E$.
As another example, consider a lemniscate
\[\sigma:=\{z:\  |T_N(z)|=1\},\]
where $T_N$ is an algebraic  polynomial of degree $N$,
in which case\footnote{To be precise, $\sigma$ consists of
Jordan curves only if $T_N'\ne 0$ on $\sigma$,
but the formula given for $\omega_E$ is true
without this assumption (excluding double points
where the density can be considered
to be 0).}
\[\omega_\sigma(z)=\frac{|T_N'(z)|}{2\pi N}.\]

If $E$ has only one component, then the equilibrium measure is
closely related to
the conformal map of its unbounded domain.
In fact, let
$E$ be a smooth Jordan curve (homeomorphic image of a circle)
or arc (homeomorphic image of a segment), and $\F$ a conformal
map from the exterior of $E$ onto the exterior of the unit circle
that leaves infinity invariant.
This $\F$ can be extended to $E$ as a continuously
differentiable function
(with the exception of the endpoints of $E$
when $E$ is a Jordan arc).
If $E$ is a Jordan curve, then simply
\[\omega_E(z)=\frac{1}{2\pi}|\F'(z)|.\]
If, however, $E$ is a Jordan arc, then it has two sides,
say positive and negative sides, and every point $z\in E$
different from the endpoints of $E$ is considered to belong
to both sides,  where they represent different points
$z_\pm$ (with different $\F$-images). In this case
\[\omega_E(z)=\frac{1}{2\pi}(|\F'(z_+)|+|\F'(z_-)|).\]
For example,  if $E$ is the arc of the unit circle that runs from
$\ee^{-\ii\beta}$ to $\ee^{\ii\beta}$ counterclockwise, then
\begin{equation}
\omega_{E}(\ee^{\ii t})=
\frac{1}{2\pi}\frac{\cos t/2}{\sqrt{\sin^{2}\beta/2-\sin^{2}t/2}},
\qquad t\in (-\beta,\beta).
\label{omegaform1}\end{equation}

\section{The general Bernstein inequality}
\subsection{Trigonometric polynomials}
The form (\ref{omegaform1}) indicates the nature of the factor in
Videnskii's inequality
(\ref{ineq110}): if $[\beta,\beta]\subset (-\pi,\pi)$,
then one must consider the arc
\[\G:=\{\ee^{\ii t}:\  t\in [-\beta,\beta]\}\]
on the unit circle, and the Videnskii factor
at a point $\theta\in (-\beta,\beta)$
is $2\pi$ times $\omega_\G(\ee^{\ii\theta})$.
It turns out that this is
true in general as is shown by  the following result
of A. Lukashov from \cite{Lukashov}
(see also \cite{TotikVidgen}).
For  a $2\pi$-periodic closed set $E\subset \R$ let
\begin{equation} \G_E:=\{\ee^{\ii t}:\  t\in E\}
\label{eg}\end{equation}
be its image when we identify $\R/({\rm mod}\ 2\pi)$
with the unit circle.
Then, for any trigonometric polynomial $T_n$ of degree at most
$n=1,2,\ldots$, we have (considering the
 one-dimensional interior ${\rm Int}(E)$ of $E$)
\begin{equation}
|T_n'(\theta)|\le n2\pi\omega_{\G_E }(\ee^{\ii\theta })\|T_n\|_E,
\qquad \theta\in {\rm Int}(E),
\label{vid1}\end{equation}
where $\omega_{\G_E}$ denotes the equilibrium density  of $\G_E$.

The result is sharp (see \cite{TotikVidgen}): if  $\theta\in E$ is
 an interior point of $E$,
then there are  trigonometric polynomials $T_n\not\equiv 0$ of
degree at most $n=1,2,\ldots$ such that
\begin{equation}
|T_n'(\theta)|\ge (1-o(1))n2\pi\omega_{\G_E }(\ee^{\ii\theta })
\|T_n\|_E,\label{vid12}\end{equation}
where $o(1)$ tends to 0 as $n\to\infty$.

\subsection{Algebraic polynomials on the real line}
The algebraic version (proved in M. Baran's paper \cite{Baran}
and independently in \cite{TotikActa}) reads as follows.
If $E\subset \R$ is a compact set,
then for algebraic polynomials $P_n$ of degree at most
 $n=1,2,\ldots$, we have
\begin{equation} |P_n'(x)|\le n{\pi\omega_E(x)}
\|P_n\|_E,\qquad x\in {\rm Int}(E).
\label{bernstein}\end{equation}
This is sharp again:
if $x_0\in {\rm Int}(E)$ is arbitrary, then there are
polynomials $P_n$ of degree at most $n=1,2,\ldots$ such that
\[ |P_n'(x_0)|\ge (1-o(1))n{\pi\omega_E(x_0)}
\|P_n\|_E.\]

\noindent We mention that (\ref{bernstein})
can also be deduced from (\ref{vid1})
via a suitable  linear transformation and
the substitution $x=\cos t$.

Note that in the special case $E=[-1,1]$
the inequality (\ref{bernstein}) gives back the
original Bernstein inequality (\ref{b1})
because $\omega_{[-1,1]}(x)=1/\pi\sqrt{1-x^2}$.

Actually, for real polynomials more than (\ref{bernstein})
is true (see \cite{TotikVidenskii}):
\begin{equation} \left(\frac{P_n'(x)}{\pi\omega_E(x)}\right)^2
+n^2P_n(x)^2\le
n^2\|P_n\|_E^2, \qquad x\in {\rm Int}(E),
\label{bernszeg}\end{equation}
which is the analogue of the beautiful inequality
\begin{equation} \left(P_n'(x)\sqrt{1-x^2}\right)^2
+n^2P_n(x)^2\le
n^2\|P_n\|_{[-1,1]}^2, \quad x\in [-1,1],\label{szeg}\end{equation}
of G. Szeg\H{o} \cite{Szego} and  G. Schaake
and J. G. van der Corput \cite{Schaake}.

\subsection{Algebraic polynomials on a circular set}
The complete analogue of (\ref{bernstein}) is known
for closed subsets $E$ of the unit
circle, see \cite{NagyTotik2}.
Indeed, if $E$ is such a set and $z\in E$ is an inner point of $E$
(i.e. an inner point of a subarc of $E$), then
for algebraic polynomials $P_n$ of degree at most $n=1,2,\ldots$
we have
\begin{equation}
|P_n'(z)|\le \frac{n}{2}\bigl(1+ 2\pi \omega_E(z)\bigr)\|P_n\|_E.
\label{balg1}\end{equation}
Furthermore, this is sharp: for an inner point  $z\in E$
there are polynomials $P_n\not\equiv0$ of degree
$n=1,2,\ldots$ for which
\[
|P_n'(z)|\ge
(1-o(1))\frac{n}{2}\bigl(1+ 2\pi \omega_E(z)\bigr)\|P_n\|_E.
\]

We shall discuss later why there is a difference in
the Bernstein factors in
(\ref{bernstein}) and (\ref{balg1}).

\section{The general Markov inequality}
\subsection{Intervals on the real line}
In the sense of the preceding section, what is the form
of the  Markov inequality (\ref{m1}) on more general
sets than an interval, say on a set consisting of finitely
many intervals?
Let $E=\cup_{i=1}^m [a_{2i-1},a_{2i}]$,
$a_1<a_2<\cdots<a_{2m}$, be such a set.
When we consider the analogue of the Markov
inequality for $E$, we actually  have to talk about a Markov-type
 local inequality  around every endpoint of $E$.
Indeed, away from the endpoints
(\ref{bernstein}) is true, therefore there
the derivative is bounded by a constant
times $n$ times the norm of the polynomial,
so  the $n^2$ factor is needed only close to
the endpoints (as in the single interval case).
It is also clear that different endpoints play different roles.
So let $a_j$ be an endpoint of $E$, and let
$E_j$ be the part of $E$ that lies closer to $a_j$ than to any
other endpoint.
Let $M_j$ be the smallest constant for which
\begin{equation}
\|P_n'\|_{E_j}\le (1+o(1))M_jn^2\|P_n\|_E,
\qquad {\rm deg}(P_n)\le n,\ n=1,2,\ldots,
\label{16a}\end{equation}
holds, where $o(1)$ tends to 0 (uniformly
in the polynomials $P_n$) as $n$ tends to infinity.
This $M_j$
depends on what endpoint $a_j$ we are considering, and it gives the
asymptotically best constant in the corresponding
local Markov inequality.
Its value can be expressed in terms of the equilibrium density
$\omega_E$. Indeed it is clear from
(\ref{omegaform}) that at $a_j$
the limit
\begin{equation}
\Omega_j:=
\lim_{t\to a_j,\ t\in E}\sqrt{|t-a_j|}\omega_E(t)
=\frac{1}{\pi}
\frac{\prod_{i=1}^{m-1}|a_j-\xi_i|}
{\sqrt{\prod_{i\not=j}|a_j-a_i|}}
\label{oj}\end{equation}
exists, where $\xi_i$ are the numbers from (\ref{7}).
For example, if $E=[-1,1]$, $a_1=-1$, $a_2=1$,
then $\Omega_{1,2}=1/\pi\sqrt 2$.
With this $\Omega_j$ the asymptotic Markov factors $M_j$
can be expressed (see \cite{TotikActa}) as
\begin{equation}
M_j=2\pi^2\Omega_j^2,\qquad j=1,\ldots,2m.
\label{16ab}\end{equation}

From here the global Markov inequality easily follows:
\begin{equation}
\|P_n'\|_E
\le (1+o(1))
n^2\left(\max_{1\le j\le 2m}2\pi^2\Omega_j^2\right) \|P_n\|_E,
\qquad {\rm deg}(P_n)\le n.\label{globalmarkov}\end{equation}
On the right  the $o(1)$ term tends to 0 uniformly
in the polynomials $P_n$ as
$n\to\infty$, and, in general, this term cannot be
dropped.
\bigskip

While the inequalities (\ref{vid1}) and (\ref{bernstein})
give the best
possible results for all $n$
(in both for each given $n$ the equality is attained
at some points), the estimates in (\ref{16a}) and
(\ref{globalmarkov}) are sharp only in an asymptotic sense because
of the term $(1+o(1))$ on the right.
Here $o(1)$ tends to zero  independently of the polynomials $P_n$
as $n\to\infty$, and the given inequality may not be true without
the $(1+o(1))$ factor.
This will be true in all subsequent results containing that factor.

If we call the $n$-th Markov constant the smallest number
$L_n=L_{n,E}$ for which
\begin{equation}
\|P_n'\|_E\le n^2L_n \|P_n\|_E
\label{globalmarkov1}\end{equation}
is true for all polynomials $P_n$ of degree at most $n$, then
the determination of $L_n$ seems to be a very difficult problem.
\begin{problem} For a given set $E$ consisting of finitely many
intervals and for a given degree $n$ find the $n$-th
Markov constant $L_n$.
\end{problem}

Analogous questions can be raised in connection with
all subsequent results that contain the $(1+o(1))$ factor
on the right. We shall not mention those problems separately.

\subsection{Markov's inequality on a system of arcs on a circle}
Suppose now that $E$ consists of finitely many circular arcs, say
\[
E=\bigcup_{k=1}^m{\{\ee^{\ii t}:\  t\in
[\alpha_{2k-1},\alpha_{2k}]\}},
\]
where $-\pi\le \alpha_1<\alpha_2<\cdots<\alpha_{2m}<\pi$.
In this case an explicit form
similar to the one in (\ref{omegaform}) is known for the
equilibrium measure (see e.g., \cite{NK2}).
We define for an endpoint $A_j=\ee^{\ii\alpha_j}$ of
a subarc of $E$ the quantity $\Omega_j$ as
\begin{equation}
\Omega_j:=\lim_{z\to A_j, \ z\in E}\sqrt{|z-A_j|}\omega_E(z).
\label{fgh}\end{equation}
With this, we have the analogue of (\ref{16a})--(\ref{16ab})
with sharp constant:
\begin{equation}
\|P_n'\|_{E_j}\le (1+o(1))
2\pi^2\Omega_j^2n^2\|P_n\|_E,
\qquad {\rm deg}(P_n)\le n,\ n=1,2,\ldots,
\label{16a*}\end{equation}
where $E_j$ is the part of $E$ that lies closer to $A_j$
than to any other endpoint in $E$.

\subsection{Markov's inequality for trigonometric polynomials}
We have already mentioned Videnskii's inequality (\ref{vid1}).
However, if
one considers derivatives of trigonometric polynomials
on an interval (or system of intervals) shorter than $2\pi$,
then a Markov-type estimate also emerges
since the factor in (\ref{vid1}) blows up around the endpoints.
Already
the original paper \cite{Videnskii1} of Videnskii contained that if
$T_n$ is a trigonometric polynomial of degree at most $n$
and $0<\beta<\pi$, then
\[
\|T_n'\|_{[-\beta,\beta]}\le
(1+o(1))n^22\cot(\beta/2)\|T_n\|_{[-\beta,\beta]}.\]

When
\[
E=\bigcup_{k=1}^m[\alpha_{2k-1},\alpha_{2k}],\qquad
-\pi\le \alpha_1<\alpha_2<\cdots<\alpha_{2m}<\pi,\]
we should again consider the set
\[\G_E=\{\ee^{\ii t}:\  t\in E\}\]
(see (\ref{eg})) and the corresponding expressions
\begin{equation}
\Omega_j^{\G_E}=
\lim_{z\to A_j, \ z\in \G_E}\sqrt{|z-A_j|}\omega_{\G_E}(z),
\qquad A_j=\ee^{\ii\alpha_j},\label{fgh1}\end{equation}
from (\ref{fgh}).
Now if $E_j$ is the part of $E$ that is closer to
the endpoint $\alpha_j$ than to any other of the endpoint
of a subinterval of $E$, then (see \cite{NK2})
\begin{equation}
\|T_n'\|_{E_j}\le
(1+o(1))n^2 8\pi^2(\Omega_j^{\G_E})^2\|T_n\|_E,
\label{sdr}\end{equation}
and here the constant on the right is sharp.

Note that in this estimate
$\pi^2(\Omega_j^{\G_E})^2$ is multiplied by 8 and
not by 2 as in the polynomials cases up to now.

\section{Jordan curves and arcs}
Let $C_1=\{z:\  |z|=1\}$ be the unit circle.
If $P_n$ an algebraic polynomial of degree at most $n$,
then  $P_n(\ee^{\ii t})$ is
a trigonometric polynomial of degree at most $n$,
so by Bernstein's inequality
(\ref{b0}), we have
\[\left|\frac{dP_n(\ee^{\ii t})}{dt}\right|\le n\max |P_n|.\]
The left-hand side is
$|P_n'(\ee^{\ii t})|\ii \ee^{\ii t}|=|P_n'(\ee^{\ii t})|$,
and we obtain
\begin{equation}
|P_n'(z)|\le n\|P_n\|_{C_1},\qquad z\in C_1.
\label{riesz}\end{equation}
This inequality is due to M. Riesz \cite{Riesz}
(although it can be easily derived from (\ref{b0}),
remember that (\ref{b0}) was originally
given with a factor $2n$ on the right, so \cite{Riesz} contains
the first correct proof of (\ref{riesz})).

\subsection{Unions of Jordan curves}
Riesz' inequality was extended to Jordan curves
and families of Jordan curves in \cite{NagyTotik}:
if $E$ is a finite union of
disjoint $C^2$-smooth Jordan curves (homeomorphic
images of circles) lying exterior to each other,
then for polynomials $P_n$ of degree at most $n=1,2,\ldots$
we have
\begin{equation}
|P_n'(z)|\le (1+o(1))n2\pi \omega_E(z)\|P_n\|_E,
\qquad z\in E.\label{gRiesz}\end{equation}
Furthermore, (\ref{gRiesz}) is best possible:
if $z_0\in E$, then there are polynomials
$P_n\not\equiv 0$ of degree at most $n=1,2,\ldots$ for which
\[|P_n'(z_0)|\ge (1-o(1))n2\pi \omega_E(z_0)\|P_n\|_E.\]

Note that if $E$ is the unit circle, then $\omega_E\equiv 1/2\pi$,
so (\ref{gRiesz}) gives back the
original inequality (\ref{riesz}) of M. Riesz modulo the
$(1+o(1))$ factor which, in general, cannot be dropped in the
Jordan curve case.

If $E$ is the union of $C^2$-smooth Jordan curves,
then (\ref{gRiesz}) implies the Markov-type norm inequality
\begin{equation}
\|P_n'\|_E\le (1+o(1))n
\Bigl(\max_{z\in E}\omega_E(z)\Bigr)\|P_n\|_E,
\label{gRieszM}\end{equation}
which is sharp again in the sense that on the
right-hand side no smaller constant
can be written than $\max\omega_E(z)$.

For the inequality (\ref{gRiesz}) at a given point $z\in E$
one does not need the $C^2$-smoothness of $E$,
it is sufficient that $E$ is $C^2$-smooth in a neighborhood of $z$.
Hence, if $E$ is the union of piecewise $C^2$-smooth
Jordan curves, then
(\ref{gRiesz}) holds at any point of $E$
which is not a corner point, i.e., where
smooth subarcs of $E$ meet.
At corner points the situation is different: if
two subarcs of $E$ meet at $z_0$ at an external
angle $2\pi\alpha$,
$0<\alpha<1$, then
\[
|P_n'(z_0)|\le Cn^{\alpha}\|P_n\|_E,\qquad {\rm deg}(P_n)\le n,
\]
see \cite{Szego0}, and here the order $n^{\alpha}$
is best possible (explaining
also why in Bernstein's inequality the order is
$n$ while in Markov's it is $n^2$).
\begin{problem} Determine the smallest $C$ for which
\[
|P_n'(z_0)|\le (1+o(1))Cn^{\alpha}\|P_n\|_E,
\qquad {\rm deg}(P_n)\le n.\]
\end{problem}

\noindent The solution of this problem would be
interesting even in such a simple case when $E$
is the unit square.

By the maximum modulus theorem both (\ref{gRiesz}) and
(\ref{gRieszM}) hold true
if $E$ is the union of the (closed) domains enclosed by
finitely many $C^2$ Jordan
curves (in which case the equilibrium measure $\mu_E$ is
supported on the boundary of $E$,
and $\omega_E$ denotes the density of $\mu_E$ with respect
to the arclength measure on that boundary).
Actually, (\ref{gRiesz}) is true under much more general
assumptions on the compact set $E$. It is sufficient
that $E$ coincides with the closure of its interior, and
its boundary is a $C^2$-smooth Jordan
arc in a neighborhood of the point $z$ where we consider
(\ref{gRiesz}).
That this is {\it not true}
when the closure assumption is not satisfied is shown
by (\ref{bernstein})
(note that there we have $\pi\omega_E(z)$ while in (\ref{gRiesz})
the correct factor is $2\pi\omega_E(z)$) and
even more dramatically by the
Jordan arc case to be discussed below.

\subsection{Bernstein's inequality on a Jordan arc}
The preceding results satisfactorily answer the form
of the Bernstein and
Markov inequalities on unions of smooth Jordan curves.
It has turned out that
the case of Jordan arcs is different and much more difficult,
and actually
we have the precise results only for one Jordan arc.
To explain why arcs are different than curves one needs
to say a few words about the
proof of (\ref{gRiesz}).
Using inverse images under polynomial maps one can
deduce (\ref{gRiesz}) from (\ref{riesz}) for lemniscates, i.e.,
sets of the form $\sigma=\{z:\  |T_N(z)|=1\}$,
where $T_N$ is a polynomial of fixed degree
(this deduction is by far not trivial,
but possible using the so-called polynomial inverse image method,
see
\cite{TotikSA})). Note that a lemniscate may
have several components, so the splitting of the underlying
domain occurs at this stage of the proof.
Now smooth Jordan curves, and actually families of Jordan curves,
can be approximated from inside
and from outside by lemniscates that touch
the set at a given point (this is done via the
sharp version of Hilbert's lemniscate theorem,
see \cite{NagyTotik}),
and that allows one to deduce (\ref{gRiesz})
in its full generality from its validity on
lemniscates $\sigma$.
Since arcs do not have interior domains, that is not
possible for arcs, and, as we shall see,
the form of the corresponding result is indeed different.

In the general inequalities we have considered so far, always
the equilibrium density $\omega_E$ gave the (asymptotically)
best Bernstein-factors, and the Markov-factors have
also been expressed in terms of them. In some sense
this was accidental, it was due to either a symmetry
(when $E\subset \R$) or to an absolute lack
of symmetry (when $E$ was a Jordan curve for which
the two sides of $E$, the exterior and interior sides of
$E$, play absolutely different roles). This is no longer
the case when we consider Jordan arcs,
for which the Bernstein factors
are not expressible via the equilibrium density.

So let $E$ be a Jordan arc, i.e., a homeomorphic image of
a segment. We assume $C^{2+\alpha}$
smoothness of $E$ with some $\alpha>0$.
As has already been discussed in Section
\ref{sectpot}, $E$  has two sides, and every point $z\in E$
different from the endpoints of $E$ gives rise
to two different points
$z_\pm$ on the two sides. With these,
$\omega_E(z)=(|\F'(z_+)|+|\F'(z_-)|)/2\pi$,
where $\F$ is a conformal
map from the exterior of $E$ onto the exterior of the unit circle
leaving infinity invariant.

Now the Bernstein  inequality on $E$ for algebraic polynomials
takes the form (for $z\in E$
being different from the two endpoints of $E$)
\begin{equation}
|P_n'(z)|\le (1+o(1))n
\max \bigl(|\F'(z_+)|,|\F'(z_-)|\bigr)\|P_n\|_E
\label{NK}\end{equation}
 (see \cite{NK}
for analytic arcs and \cite{NKT}
for the general case).
This is best possible:
one cannot write anything smaller than
$\max (|\F'(z_+)|,|\F'(z_+)|)$
on the right. As for the $o(1)$ term in (\ref{NK}), it
 may depend on the position of $z$ inside $E$, but
it is uniform in $z$ on any  closed subarc of $E$
that does not contain the endpoints of
$E$ and, as before,  and it is uniform in $P_n$, $n=1,2,\ldots$.

\begin{figure}
\begin{center}
\includegraphics[width=.3\textwidth]{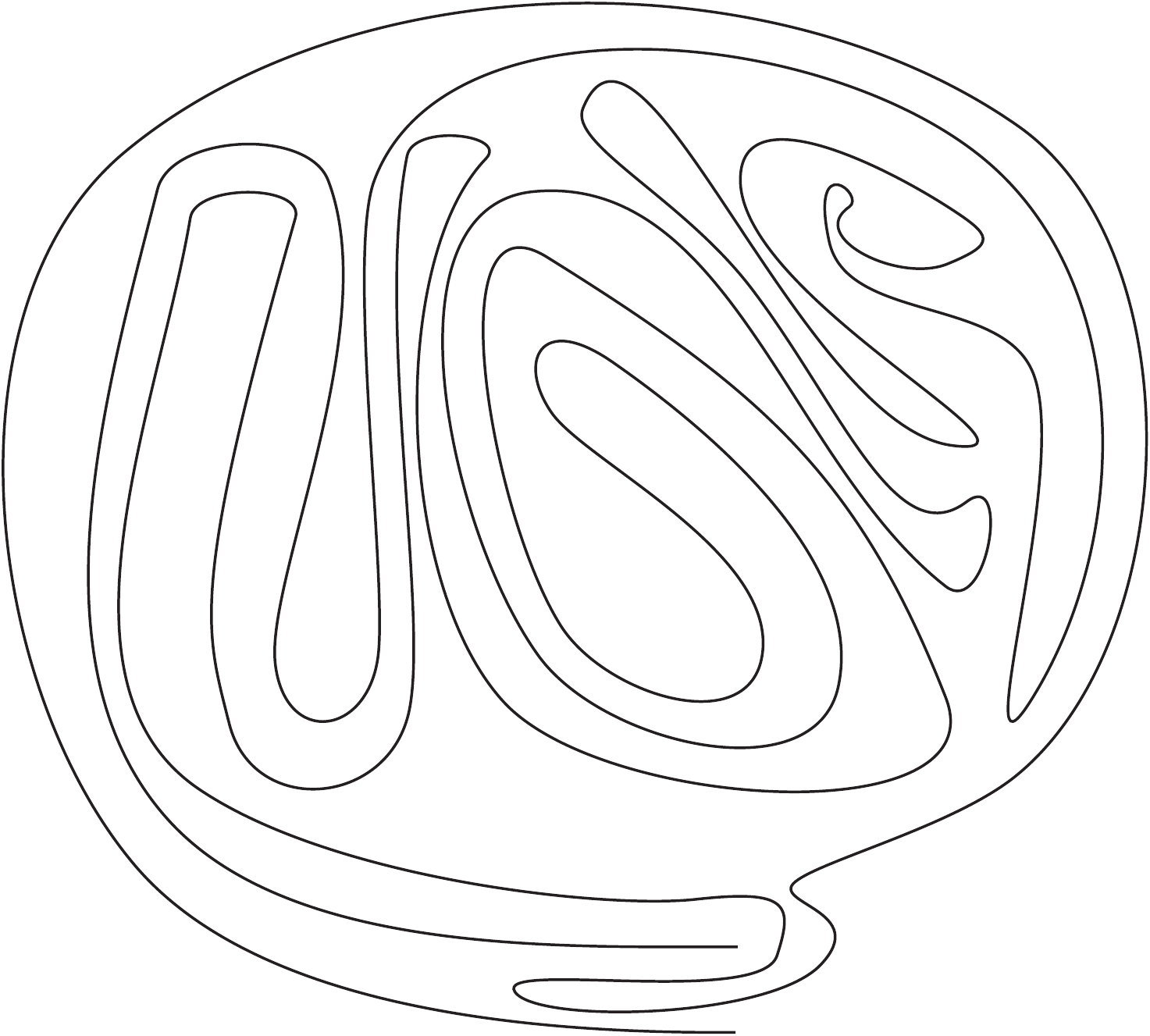}
\end{center}
\caption{A ``wild" Jordan arc\label{fig2}}
\end{figure}

To appreciate the strength of (\ref{NK})
(or that of (\ref{gRiesz}))
 let us mention that  the (smooth) Jordan arc
in it can be arbitrary, and
a general (smooth) Jordan arc can be pretty complicated,
see for example,
Figure \ref{fig2}.
\bigskip

\begin{problem} Find the analogue of (\ref{NK})
for $E$ consisting of more than one
(smooth) Jordan arc or when $E$ is the union of
 Jordan curves and arcs.
\end{problem}

We believe that the answer to this problem will be the following.
There is a possibly multivalent analytic function
$\Psi$ in the unbounded component
$\Omega$ of $\C\setminus E$ that maps $\Omega$
onto the exterior of the unit circle
locally conformally. While this $\Psi$ is multivalent,
its absolute value
$|\Psi|$ is single-valued, and $g_\Omega(z)=\log|\Psi(z)|$
is actually the Green's function
of $\Omega$ with pole at infinity
(there are other definitions of the Green's function,
one should take any of them).
Now (in the one component case when $\Psi$ is just the
conformal map $\F$ that was considered before)
the moduli $|\F'_\pm(z)|$ in (\ref{NK})
are precisely the normal derivatives
$\partial g_\Omega(z)/\partial {\bf n}_\pm$ of $g_\Omega$
in the direction of
the two normals to $E$ at $z$, hence (\ref{NK}) can be written as
\begin{equation}
|P_n'(z)|\le (1+o(1))n
\max \left(\frac{\partial g_\Omega(z)}{\partial {\bf n}_+},
\frac{\partial g_\Omega(z)}{\partial {\bf n}_-}\right)\|P_n\|_E,
\label{NK11}\end{equation}
and it is expected that this form remains true not
just when $E$ is a single Jordan arc, but also when
 $E$ is the union of smooth Jordan arcs and curves
(if $z$ belongs to a Jordan curve,
then the normal derivative in the direction of
the inner domain is considered 0).

The conjecture just explained is true in two special cases:
 when $E$ is a union of
real intervals and when $E$ is the union of finitely many arcs
 on the unit circle.
In fact, both in (\ref{bernstein}), resp. (\ref{balg1}),
that cover these cases,
the Bernstein factors $\pi\omega_E(x)$, resp.
$(1+ 2\pi \omega_E(z))/2$, are precisely the maximum of the
normal derivatives (see \cite{NagyTotik2}).
\medskip

If $E$ is a piecewise smooth Jordan curve which
may have ``corners", then (\ref{NK}) still holds
for points where $z$ is smooth.

\begin{problem} Find the analogue of (\ref{NK})
for a piecewise smooth Jordan arc $E$ at
corner points.
\end{problem}
If at a corner point the two connecting subarcs form
complementary angles $\alpha 2\pi$ and $(1-\alpha)2\pi$,
$0\le \alpha\le 1/2$, then
\[
|P_n'(z)|\le
(1+o(1))Cn^{1-\alpha}\|P_n\|_E,\qquad {\rm deg}(P_n)\le n,
\]
with some constant $C$, and the problem is
to determine the smallest $C$.
This is not known even in such simple cases when $E$
is the union of two perpendicular segments of
equal length.

\subsection{Markov's inequality on a Jordan arc}
As for Markov's inequality, let the endpoints of the
Jordan arc $E$ be the points $A$ and $B$.
Consider e.g., the endpoint $A$, and  let $\tilde E$
be the part of $E$ that is closer to $z$ than to the other
endpoint of $E$.
It turns out that as $z\to A$ the density $\omega_E(z)$
behaves like $1/\sqrt{|z-A|}$, and actually the  limit
\[
\Omega_A:=\lim_{z\to A,\ z\in E} \sqrt{|z-A|}\omega_E(z)
\]
exists. With it we have the
Markov inequality around $A$ (see \cite{TotikSb}):
\begin{equation}
\|P_n'\|_{\tilde E}\le
(1+o(1))n^2 2\pi^2\Omega_A^2\|P_n\|_E,\qquad {\rm deg}(P_n)\le n,
\label{mm}\end{equation}
and this is best possible in the sense that
 one cannot write a smaller
number than $2 \pi^2\Omega_A^2$ on the right.
From here the global Markov inequality
\begin{equation}
\|P_n'\|_{E}\le
 (1+o(1))n^2 2\pi^2(\max(\Omega_A,\Omega_B))^2\|P_n\|_E,
\qquad {\rm deg}(P_n)\le n,
\label{mmglobal}\end{equation}
follows immediately, and this is sharp again.

\begin{problem} Prove (\ref{mm}) when $E$ is a
union of smooth Jordan arcs.
\end{problem}
That (\ref{mm}) should be the correct form also for a
system of arcs is indicated
by (\ref{16ab}) and (\ref{16a*}) which are the special
cases when $E$ is the union of finitely many
intervals or the union of finitely many arcs on the unit circle.

\bigskip
We also mention that the just discussed results in this section
are valid in a suitable form
not only for polynomials, but also for rational functions
for which the poles stay away from $E$, see \cite{NKT}.
\section{Higher derivatives}
For higher derivatives the correct form of the
Markov inequality (\ref{m1}) was given
in 1892 by V. A. Markov \cite{Markov2}, the brother of
A. A. Markov: if $k\ge 1$ is a natural number, then
\begin{equation}
\|P_n^{(k)}\|_{[-1,1]}\le
\frac{n^2(n^2-1^2)(n^2-2^2)\cdots(n^2-(k-1)^2)}
{1\cdot 3\cdots (2k-1)}
\|P_n\|_{[-1,1]}.
\label{20}\end{equation}
The  equality is attained
for the Chebyshev polynomials
$P_n(x)=\cos(n\arccos x)$.
If we write (\ref{20}) in the less precise form
\begin{equation}
\|P_n^{(k)}\|_{[-1,1]}\le
\frac{n^{2k}}{(2k-1)!!}\|P_n\|_{[-1,1]},
\label{201}\end{equation}
and compare it with
\[\|P_n^{(k)}\|_{[-1,1]}\le n^{2k}\|P_n\|_{[-1,1]}\]
which is obtained from the original Markov inequality (\ref{m1})
by iteration, then
we can see a mysterious improvement of $1/(2k-1)!!$.
It turns out that the same
improvement appears in other higher order Markov-type
inequalities, as well, but that is not
the case for Bernstein-type estimates.

\subsection{Higher order Markov inequalities}
Indeed, let $E=\cup_{i=1}^m[a_{2i-1},a_{2i}]$ be
a set consisting of finitely many intervals. Then the analogue of
of (\ref{16a})--(\ref{16ab}) for higher derivatives
 is (see \cite{TotikZ})
\begin{equation}
\|P_n^{(k)}\|_{E_j}\le
(1+o(1))\frac{2^k\pi^{2k}\Omega_j^{2k}}{(2k-1)!!}\|P_n\|_E
\label{16ak}\end{equation}
with an asymptotically sharp factor on the right.
 From here the global Markov inequality
\begin{equation}
\|P_n^{(k)}\|_{E}\le
(1+o(1))\frac{2^k\pi^{2k}(\max_j\Omega_j)^{2k}}
{(2k-1)!!}\|P_n\|_E
\end{equation}
is an easy consequence.

In a similar fashion, if $E$ is a Jordan arc as in
 (\ref{mm}) with endpoints $A$ and $B$,
then we have (see \cite{NKT}, \cite{TotikSb})
\begin{equation}
 \|P_n^{(k)}\|_{\tilde E}\le
 (1+o(1))n^{2k}\frac{2^k\pi^{2k}\Omega_A^{2k}}{(2k-1)!!}\|P_n\|_E,
\label{jkl}\end{equation}
and, as an immediate consequence,
\[
\|P_n^{(k)}\|_{E}\le
(1+o(1))n^{2k}
\frac{2^k\pi^{2k}\max(\Omega_A,\Omega_B)^{2k}}
{(2k-1)!!}\|P_n\|_E,
\]
again with the best constants
(i.e., no smaller number can be written on the right).

We do not have an explanation for the factor $1/(2k-1)!!$,
but we do know how it appears.
Consider e.g., (\ref{jkl}), and assume that the endpoint $A$ is
at the origin. Then
\[\G:=\{z:\  z^2\in E\}\]
is a Jordan arc symmetric with respect to  the origin
for which $0$ is an ``inner" point, and for it the
quantities $|\F'(0_\pm)|$ from (\ref{NK}) are the same,
and can be expressed by $\Omega_0$.
Consider $R_{2n}(z):=P_n(z^2)$. For $k\ge 2$
the term  $P_n^{(k)}(z^2)$ appears in the $2k$-th
derivative of $R_{2n}(z)=P_n(z^2)$
if we use Fa\'a di Bruno's formula for the $2k$-th
derivative of composite functions, and
 $1/(2k-1)!!$ appears as the coefficient of that term
(when everything is evaluated at $z=0$).
Now an application of
(\ref{NKk}) below with $2k$ instead of $k$
for $R_{2n}$ at $z=0$ yields the bound given in (\ref{jkl})
(at least at the endpoint 0).

\subsection{Higher order Bernstein inequalities}
When we consider Bernstein-type estimates the situation is
different, no improvement factor like
$1/(2k-1)!!$ appears.
Indeed,  the higher derivative form of (\ref{bernstein})
and  (\ref{NK}) are
\begin{equation}
|P_n^{(k)}(x)|\le (1+o(1))n^k({\pi\omega_E(x)})^k
\|P_n\|_E,\qquad x\in {\rm Int}(E),
\label{bernsteink}\end{equation}
(when $E\subset \R$), and
\begin{equation}
|P_n^{(k)}(z)|
\le (1+o(1))n^k
\max \bigl(|\F'(z_+)|,|\F'(z_-)|\bigr)^k\|P_n\|_E
\label{NKk}\end{equation}
(when $E$ is a Jordan arc), which are best possible.
So in these cases the best results are obtained from the
estimate on the first derivative  by taking formal powers,
and there is no improvement of the sort
$1/(2k-1)!!$  as opposed to the above-discussed
Markov inequalities.

In a similar manner, if $E$ consists of a finite number of
smooth Jordan curves,
then the Riesz inequalities (\ref{gRiesz}) and (\ref{gRieszM})
for higher
derivatives take the best possible forms
\begin{equation}
|P_n^{(k)}(z)|\le
(1+o(1))n^k(2\pi \omega_E(z))^k\|P_n\|_E,
\qquad z\in E,\label{gRieszk}\end{equation}
and
\[
\|P_n^{(k)}\|_E\le
(1+o(1))n^k\Bigl(\bigl(\max_{z\in E}\omega_E(z)\bigr)\Bigr)^k
\|P_n\|_E,
\]
so there is no improvement again compared to straight iterations.

While (\ref{bernsteink})--(\ref{gRieszk}) seem to appear as
iterations of the $k=1$ case,
no straightforward iteration is possible.
However, the proofs still use the $k=1$ case inductively
 in combination with a localization technique using so-called
fast decreasing polynomials.

We close this section by stating the higher order analogue of
(\ref{vid1}), i.e., the
higher order Bernstein inequality for trigonometric polynomials:
if $E\subset \R$ is a $2\pi$-periodic
closed set, then
\begin{equation}
 |T_n^{(k)}(\theta)|\le
 n^k\Bigl(2\pi\omega_{\G_E }(\ee^{\ii\theta })\Bigr)^k\|T_n\|_E,
\qquad \theta\in {\rm Int}(E),
\label{vid1k}\end{equation}
where $\omega_{\G_E}$ denotes the equilibrium density  of the
set (\ref{eg}).
As before, (\ref{vid1k}) is sharp in the sense that
no smaller factor than
$(2\pi\omega_{\G_E }(\ee^{\ii\theta }))^k$ can be written
on the right.

The inequality (\ref{NKk}) appears in \cite{NKT},
(\ref{vid1k}) in \cite{NK1},
and the proof of (\ref{bernsteink}) was given in
Appendix 2 of \cite{TotikA}.
While (\ref{gRieszk}) has not been recorded before,
it can be deduced from the $k=1$ case using the machinery of
\cite{NK2} or \cite[Appendix 2]{TotikA} .
\bigskip

The higher order versions of (\ref{balg1}), (\ref{16a*})
and (\ref{sdr})
 are also known and follow the above pattern ($1/(2k-1)!!$
improvement in the Markov case and no improvement
in the Bernstein case).
We refer the reader to \cite{NK2}.

\section{\texorpdfstring{$L^2$}{L2}-Markov inequalities}
The $L^p$ version of the preceding results is much less known.
Here we shall consider
only a few results mostly related to the case $p=2$.

Let $\nu_\kappa$ denote the smallest positive zero of the
Bessel function $J_{(\kappa-1)/2}$ of the first kind
(see e.g., \cite{Watson}).
It was proved in  \cite{Aptekarev}
by  A.~I.~Aptekarev, A.~Draux, V.~A.~Kalyagin and D.~N.~Tulyakov
 that for polynomials $P_n$ of degree at most $n=1,2,\ldots$
\begin{equation}
\left(\int_{-1}^1|P_n'(x)|^2dx\right)^{1/2}\le
(1+o(1))n^2\frac{1}{2\nu_0}
\left(\int_{-1}^1|P_n(x)|^2dx\right)^{1/2}.
\label{Apt0}\end{equation}
Furthermore, on the right $1/2\nu_0$ is the
smallest possible constant.

If $E=\cup_{i=1}^m[a_{2i-1},a_{2i}]$ is the union of $m$
intervals, then
the extension of (\ref{Apt0}) to $E$ reads  as
(see \cite{TotikL2})
\begin{equation}
\left(\int_E|P_n'(x)|^2dx\right)^{1/2}
\le (1+o(1))
n^2\bigl(\max_j\pi^2\Omega_j^2\bigr)\frac{1}{\nu_0}
\left(\int_E|P_n(x)|^2dx\right)^{1/2},
\label{Apt01}\end{equation}
where $\Omega_j$ are the quantities defined in (\ref{oj}).
Furthermore, this estimate is
sharp, no smaller constant can be written on the right.

Recall that if $E=[-1,1]$, $a_1=-1$, $a_2=1$,
then $\Omega_{1,2}=1/\pi \sqrt 2$, so
(\ref{Apt01}) reduces to (\ref{Apt0}).

More generally, let $w(x)=(1+x)^\alpha(1-x)^\beta$,
$\alpha,\beta>-1$,
be a Jacobi weight. Then
the sharp $L^2$-Markov inequality with this weight is
\begin{multline} 
\left(\int_{-1}^1|P_n'(x)|^2w(x)dx\right)^{1/2}
\\ \le (1+o(1))
n^2\frac{1}{2\nu_{{\rm min}(\alpha,\beta)}}
\left(\int_{-1}^1|P_n(x)|^2w(x)dx\right)^{1/2}
\label{Apt}
\end{multline}
(see \cite{Aptekarev} for $|\alpha-\beta|\le 4$
and \cite{TotikL2} for the other cases).

The analogue of this for several intervals is as follows.
Let 
$E=\cup_{i=1}^m[a_{2i-1},a_{2i}]$ be  the union of $m$ 
intervals and
\[w(t)=h(t)\prod_{i=1}^{2m} |t-a_i|^{\alpha_i},\]
$\alpha_i>-1$, a generalized Jacobi weight,
where $h$ is a positive continuous function on $E$.
Then (see \cite{TotikL2})
\begin{equation}
\|P_n'\|_{L^2(w)}
\le (1+o(1))n^2M(E,w)\|P_n\|_{L^2(w)},
\qquad {\rm deg}(P_n)\le n,
\label{L1}\end{equation}
where the smallest possible constant $M(E,w)$ is
\begin{equation}
M(K,w)=\max_{1\le j\le 2m} \frac{\pi^2 \Omega_j^2}{\nu_{\alpha_j}}.
\label{L2}\end{equation}

\begin{problem}
Find the precise form of these inequalities in other $L^p$,
$1\le p<\infty$, norms.
\end{problem}

The Bernstein-type version of (\ref{Apt0})/(\ref{Apt}) was found by
A.~Guessab and G.~V.~Milovanovic
\cite{Milovanovic} much earlier and actually in a stronger form:
if 
$w(x)=(1+x)^\alpha (1-x)^\beta$, $\alpha,\beta>-1$, is a Jacobi weight, then 
\begin{multline} 
\left(\int_{-1}^1|\sqrt{1-x^2}P_n'(x)|^2w(x)dx\right)^{1/2}
\\\le 
\sqrt{n(n+1+\alpha+\beta)}
\left(\int_{-1}^1|P_n(x)|^2w(x)dx\right)^{1/2},
\label{Grad}
\end{multline}
with equality for the corresponding Jacobi polynomial of degree $n$.
Remarkably, \cite{Milovanovic}
contains also the analogue of this inequality for higher derivatives
as well with precise
constants for all $n$.

In the $\alpha=\beta=-1/2$ case this can be written in the somewhat 
less precise form
\begin{multline} 
\left(\int_{-1}^1
\left|\sqrt{1-x^2}  P_n'(x)\right|^2
\frac{1}{\sqrt{1-x^2}}
dx\right)^{1/2}
\\\le
n(1+o(1))
\left(\int_{-1}^1|P_n(x)|^2\frac{1}{\sqrt{1-x^2}}dx\right)^{1/2},
\label{Grad1}
\end{multline} 
and this form it has an extension to other $L^p$ spaces and to 
several intervals (see
\cite{NagyTookos}): let $E\subset \R$ be a compact set consisting 
of non-degenerate intervals.
Then for $1\le p<\infty$ and for algebraic polynomials $P_n$ of 
degree $n=1,2,\ldots$ we have
\[
\left(\int_E\left|\frac{P_n'(x)}{\pi\omega_E(x)}\right|^p
\omega_E(x)dx\right)^{1/p}
\le n(1+o(1))
\left(\int_E|P_n(x)|^p\omega_E(x)dx\right)^{1/p},
\]
and this is precise in the usual sense.
Note that if $E=[-1,1]$, then $\pi\omega_E(x)=1/\sqrt{1-x^2}$,
so in this case this inequality reduces to (\ref{Grad1})
for $p=2$.

\bigskip

\noindent {\bf Acknowledgement.} In the original version of
this paper the $\alpha=\beta=0$ 
case of (\ref{Grad}) was proposed as a problem,
and we thank
G. Milovanovic for pointing out \cite{Milovanovic}
where the solution can be found.

\bigskip

S. Kalmykov \\
School of Mathematical Sciences\\
Shanghai Jiao Tong University \\
800 Dongchuan Rd \\
Shanghai 200240, China

and

\noindent
Institute of Applied Mathematics, FEBRAS\\
7 Radio St.\\
Vladivostok, 690041, Russia\\
{\tt kalmykovsergei@sjtu.edu.cn}

\bigskip

B. Nagy\\
MTA-SZTE Analysis and Stochastics Research Group\\
Bolyai Institute\\
University of Szeged\\
Szeged\\
Aradi v. tere 1, 6720, Hungary\\
{\tt nbela@math.u-szeged.hu}

\bigskip

V. Totik\\
MTA-SZTE Analysis and Stochastics Research Group\\
Bolyai Institute\\
University of Szeged\\
Szeged\\
Aradi v. tere 1, 6720, Hungary\\
{\tt totik@math.u-szeged.hu}

\end{document}